\documentclass[Conference]{IEEEtran}

\usepackage[margin=0.75in]{geometry}

\makeatletter
\def\ps@IEEEtitlepagestyle{%
  \def\@oddfoot{\mycopyrightnotice}%
  \def\@evenfoot{}%
}
\def\mycopyrightnotice{%
  {\footnotesize 979-8-3315-5720-1/26/\$31.00~\copyright~2026 IEEE\hfill}
  \gdef\mycopyrightnotice{}
}

\usepackage{multirow}
\usepackage{multicol}
\usepackage{lipsum}
\usepackage{graphicx}
\usepackage{graphics}
\usepackage{amsmath}
\usepackage{amsbsy}
\usepackage{blindtext}
\usepackage{tcolorbox}
\usepackage{amssymb}
\usepackage{scalerel}
\usepackage{mathabx}
\usepackage{verbatim}
\usepackage{booktabs}
\usepackage{rotating,tabularx}
\usepackage{mathrsfs} 
\usepackage{hyperref}

\usepackage{adjustbox}
\usepackage{soul}
\usepackage{xcolor}
\usepackage{cite}
\usepackage{tikz}
\usepackage{color, colortbl}
\usepackage[first=0,last=9]{lcg}
\definecolor{LightGray}{rgb}{0.7,0.7,0.7}

\usepackage[wby]{callouts}
\usetikzlibrary{shapes.multipart}

\usepackage{array}
\usepackage{makecell}

\allowdisplaybreaks

\usepackage{amsthm}
\theoremstyle{definition}

\theoremstyle{remark}

\usepackage[utf8]{inputenc}
\usepackage[english]{babel}

\usepackage[caption=false,font=footnotesize]{subfig}
\usepackage[T1]{fontenc}
\usepackage{scalerel,stackengine}
\newcommand\reallywidecheck[1]{%
\savestack{\tmpbox}{\stretchto{%
  \scaleto{%
    \scalerel*[\widthof{\ensuremath{#1}}]{\kern-.6pt\bigwedge\kern-.6pt}%
    {\rule[-\textheight/2]{1ex}{\textheight}}
  }{\textheight}%
}{0.5ex}}%
\stackon[1pt]{#1}{\scalebox{-1}{\tmpbox}}%
}
\stackMath
\IEEEoverridecommandlockouts
\IEEEoverridecommandlockouts

\makeatletter
\patchcmd{\@maketitle}
{\addvspace{0.5\baselineskip}\egroup}
{\addvspace{-1\baselineskip}\egroup}
{}
{}
\makeatother

\newcommand*{\ra}{\textcolor{black}}

\usepackage{bm}

\newif\ifarxiv
\arxivtrue
\arxivfalse

\renewcommand{\baselinestretch}{.96}

\begin{document}

\title{\LARGE\bf
A Transformer-Based Mixture-of-Experts Framework for False Data Injection Attack Detection and Localization}

\author{\ra{Ruslan Abdulin,$^{\ast}$ Mohammad Rasoul Narimani$^{\ast}$
\thanks{\ra{${*}$: Department of Electrical and Computer Engineering, California State University, Northridge (CSUN). ruslan.abdulin.673@my.csun.edu, Rasoul.narimani@csun.edu. Support from NSF contract \#2523881.}%
}}}

\maketitle

\begin{abstract}
\ra{False data injection attacks (FDIAs) threaten smart-grid operation by manipulating measurement data and misleading power system state estimation. Although recent data-driven methods have shown promising performance, many existing approaches rely on a single graph-filtering mechanism and therefore struggle to adapt to diverse cyberattack patterns. This paper proposes a topology-aware encoder-only transformer mixture-of-experts framework for joint FDIA detection and localization. The proposed model integrates Laplacian positional encoding and diffusion-biased self-attention to capture global spatial dependencies, while a densely routed mixture-of-experts module composed of Chebyshev and autoregressive moving-average graph convolutional experts adaptively applies complementary graph filters based on the current grid state. The framework is evaluated on the IEEE 118- and 300-bus systems using the New York Independent System Operator load profile. Experimental results demonstrate that the proposed approach achieves up to 93.91\% detection F1 score and 84.97\% localization F1 score, outperforming existing benchmark methods while maintaining low false-alarm rates.}
\end{abstract}

\section{Introduction}
\label{Introduction}

\ra{Smart power grids are complex cyber-physical systems that rely on communication networks and software-driven energy management systems (EMSs) to maintain stable operation through power system state estimation (PSSE) algorithms. Because operational decisions depend heavily on measurement integrity, compromised sensor data can directly affect the estimated network state and overall grid stability. False data injection attacks (FDIAs) exploit this vulnerability by manipulating smart-grid measurements to mislead monitoring and control processes. Although bad data detection (BDD) mechanisms are commonly employed to identify malicious measurements, traditional residual-based approaches can be bypassed when compromised data remain consistent with the measurement model and produce sufficiently small residuals~\cite{pourshirazi2025false}. Consequently, accurate attack localization is as important as attack detection, since identifying compromised buses enables system operators to isolate affected regions and mitigate cascading impacts. Due to the severe risks posed by FDIAs, joint detection and localization remain critical research challenges in modern smart grids~\cite{boyaci2021joint}.}

\ra{Current approaches for handling FDIAs can generally be categorized into model-based and data-driven methods \cite{pourshirazi2025false}. Model-based approaches rely on system modeling, parameter estimation, and physical grid characteristics to identify inconsistencies in measurements using power flow equations and real-time observations~\cite{wang2024fdia}. However, their performance often depends on the accuracy of the underlying system model, which limits scalability and applicability \cite{kecceci2025distributed}. In contrast, data-driven methods learn to distinguish between normal and compromised operating conditions directly from historical data, offering greater flexibility and scalability for practical smart-grid applications \cite{fahim2024generalized}. Motivated by advances in machine learning (ML), a wide range of ML-based approaches have been investigated for FDIA detection and localization, including traditional ML methods, deep learning architectures, temporal models, graph neural networks, attention-based frameworks, federated learning techniques, and hybrid combinations thereof \cite{Drayer2019GraphModulation, Hasnat2020GraphSignalProcessing, boyaci2022cyberattack, li2024detection}.}

\ra{Although graph neural networks and attention-based architectures, such as the autoregressive moving average convolutional (ARMAConv) operator, Chebyshev graph convolutional (CGCN) network, and ARMAConv Encoder-only Transformer (ACEOT), have been applied to FDIA detection and localization, many existing methods rely on a single graph-filtering operator or a fixed inductive bias. As a result, their performance may degrade when attack patterns require adaptive graph-filtering behavior. For example, Chebyshev-based techniques strongly depend on the polynomial order $K_{\mathrm{Cheb}}$, where small values may miss multi-hop inconsistencies while large values can blur localization performance \cite{defferrard2016convolutional}. Moreover, practical cyber intrusions often evolve through multiple stages, including reconnaissance, attack execution, and stealth maintenance phases \cite{lian2022critical}. Since adversaries continuously adapt their attack strategies to maximize impact and evade detection \cite{yang2023multi}, FDIA handling requires models capable of dynamically adjusting their analytical behavior.}

\ra{Mixture of Experts (MoE) models address this challenge by combining specialized submodels, called experts, whose contributions are selected according to the current grid state. Consequently, as attack dynamics change, the most suitable expert can make the dominant contribution to FDIA detection and localization. MoE frameworks can employ either dense (soft) or sparse (hard) routing. In dense routing, the gating mechanism assigns weights to all experts, whereas sparse routing activates only a subset of experts \cite{jacobs1991adaptive}. Although sparse routing can reduce computational cost, sparse backward updates may introduce training instability and suboptimal performance, motivating the use of denser routing when computational overhead is acceptable \cite{nie2021dense}.}

\ra{Despite the growing interest in MoE applications, particularly in large language models such as Wu Dao 2.0 and Mixtral 8x7B \cite{jiang2024mixtral, lo2025closer}, the application of multi-operator architectures to FDIA handling remains limited. To address this gap, this work proposes an Encoder-only Transformer Convolution Mixture-of-Experts (EOTConvMoE) model for joint FDIA detection and localization. The framework integrates Laplacian positional encoding and diffusion-biased self-attention to generate globally contextualized node embeddings before expert routing. We also investigate both dense and sparse routing strategies and augment the loss function using entropy regularization and router log-partition penalties to improve routing stability and balanced expert utilization. Finally, the proposed method is evaluated against common deep learning benchmarks under identical hyperparameter optimization settings on the IEEE 118- and 300-bus systems using the NYISO July 2021 load profile.}

\ra{The remainder of this paper is organized as follows: Section \ref{Problem Formulation} defines the problem, Section \ref{Proposed method} presents the proposed method and architecture, Section \ref{Experiments} presents the experimental setup and numerical results, and Section \ref{Conclusion} concludes the paper.}

\section{Problem Formulation}
\label{Problem Formulation}
\subsection{System protection and FDIA targets}

\ra{A power network state can be represented by $\bm{x}$, which contains the voltage magnitude $V_i$ and phase angle $\theta_i$ at each bus $i$. The weighted least squares estimation (WLSE) algorithm is commonly used to obtain the best-fit state estimate $\bm{\hat{x}}$ from active and reactive power injections and power flow measurements $\bm{z}$ collected by Remote Terminal Units (RTUs) and Phasor Measurement Units (PMUs) \cite{boyaci2021joint}}

\begin{equation}
    \hat{\bm{x}} = \arg\min_{\bm{x}} (\bm{z} - h(\bm{x}))^{T} \bm{R}^{-1} (\bm{z} - h(\bm{x}))
    \label{eq:wlse}
\end{equation}

\ra{where $h(\cdot)$ is the nonlinear measurement function mapping a system state to predicted measurements, $\bm{R}$ is the measurement error covariance matrix, and $\bm{z}$ contains active and reactive power injections and flow measurements, including $P_i$, $Q_i$, $P_{ij}$, and $Q_{ij}$.}

\ra{Residual-based BDD is commonly used as the first line of defense against malicious measurements in energy management systems. The residual vector $\bm{r}$ is computed by subtracting the predicted measurements $h(\bm{\hat{x}})$ from $\bm{z}$, after which the BDD module evaluates whether the residual magnitude indicates a potential FDIA. In the Largest Normalized Residual Test (LNRT)-based BDD, the largest normalized residual is compared against a threshold $\tau$ determined by the $3\sigma$ rule \cite{asefi2023anomaly}:}

\ra{\begin{equation}
    \max_{i}|\frac{r_i}{\sqrt{\Omega_{ii}}}| \ge \tau
\end{equation}
where $\Omega$ denotes the residual covariance matrix.}

\ra{An adversary aims to manipulate measurements while maintaining residuals small enough to bypass BDD. To achieve this, FDIA injects an attack vector $\bm{a}=h(\bm{x}_a)-h(\bm{x}_0)$ into healthy measurements $\bm{z}_0$, producing compromised measurements $\bm{z}_a=\bm{z}_0+\bm{a}$. As a result, the biased estimate $\bm{\hat{x}}_a$ yields a small residual $\bm{r}_a=\bm{z}_a-h(\bm{\hat{x}}_a)$, allowing the attack to remain consistent with the measurement model and evade detection \cite{liu2011false}.}

\subsection{Problem classification and constraints}

\ra{The objective of FDIA handling is to detect attacks and identify compromised buses. Accordingly, FDIA detection and localization can be formulated as a multi-label classification problem with two label types: graph-level labels for attack detection and node-level labels representing the status of each bus. For simplicity, we use binary labels, where 0 denotes a healthy network/bus, and 1 indicates an attacked graph/node.}

\ra{In practice, attackers typically have limited access to meters and devices and must carefully coordinate compromised measurements to maximize impact while minimizing exposure. As a result, FDIAs are often localized but can still cause significant economic damage \cite{zhang2018can}. To reflect these practical constraints, we assume that attacks affect regions of interconnected measurements.}

\subsection{Graph-based approach to the problem}

\ra{To construct the attack vector $\bm{a}$, an adversary must possess knowledge of the grid topology, node attributes, and electrical coupling. Consequently, FDIA handling can be formulated as a graph-based problem. For a smart grid with $n$ buses, the system is represented as a weighted, undirected graph $\mathcal{G}=(\mathcal{V},\mathcal{E},\mathbf{A}_w)$, where $\mathcal{V}$ denotes the set of buses, $\mathcal{E}$ represents the physical branches, and $\mathbf{A}_w\in\mathbb{R}^{n\times n}$ is the weighted adjacency matrix derived from the bus admittance matrix.}

\ra{Due to the irregular topology of power networks, conventional methods such as convolutional neural networks (CNNs) cannot be directly applied because they rely on Euclidean structures and consistent local neighborhoods \cite{qu2023hyperbolic}. Therefore, FDIA handling requires graph-based techniques capable of processing non-Euclidean data.}

\section{Proposed method}
\label{Proposed method}

\ra{This section presents the proposed framework for joint graph-based FDIA detection and localization. The proposed EOTConvMoE model addresses two key challenges: localization under varying graph spectral characteristics and the limited receptive scope of convolutional message-passing operators. The framework combines topology-aware attention with a densely routed MoE module that adaptively leverages CGCN and ARMAConv experts. The following subsections describe the key components and architecture of the proposed model.}

\subsection{The Mixture-of-Experts concept}

\begin{figure}
    \centering
    \includegraphics[width=0.45\textwidth]{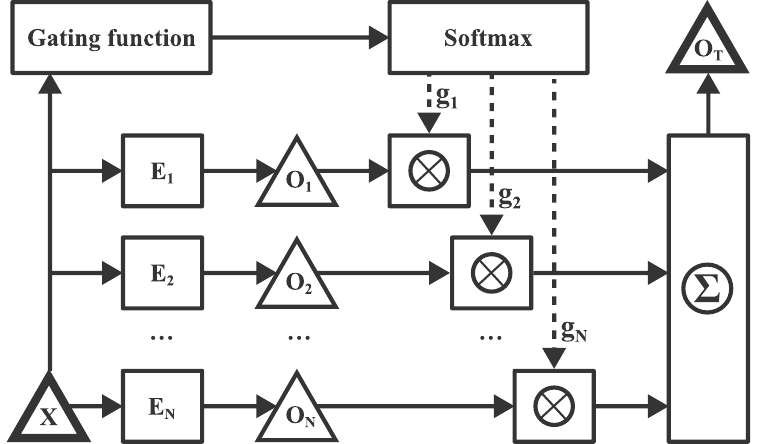}
    \caption{General MoE Architecture with the number of experts $N \ge 2$. The input is independently shared with every expert and the routing mechanism. The router produces normalized expert weights through the softmax functions defined in Eq.~\ref{eq:MOE_coef}. The final output is the weighted sum of expert outputs and the corresponding contributions defined in Eq.~\ref{eq:MOE_coef}. The resulting representation $O_T$ has the same dimensionality as $O_i$. }
    \label{fig:generalMoE}
\end{figure}

\ra{The Mixture of Experts (MoE) method is a model combination technique that uses multiple data-driven architectures, called experts, together with a gating mechanism to generate the final output. In general, an MoE model consists of two main components: (i) a finite set of learnable functions, referred to as experts, and (ii) a gating function that determines how much each expert contributes to the final prediction. Although the gating function can assign different weights to individual nodes within each expert, this study uses a graph-level router function that assigns a single scalar value $\mathrm{g_i}$ to each expert $i$.}
\ra{A common example of an MoE architecture is the Mixture of Multilayer Perceptron experts (MME), which employs a dense gating mechanism \cite{ebrahimpour2008view}. During inference, the MME architecture first processes the input through a network of MLP layers followed by a softmax activation function to produce a vector of nonnegative mixture coefficients, denoted by $\mathrm{g}$. These coefficients can be interpreted as estimated probabilities representing the contribution of each expert $i$ to the final MoE representation. The coefficients are then multiplied by the outputs of the corresponding experts to control their influence, and the weighted outputs are summed to produce the final prediction \cite{masoudnia2014mixture}. This process is illustrated in Figure \ref{fig:generalMoE} and can be formally expressed as:}

\begin{align}
\mathrm{g_i} =\frac{\exp(\ell_i)}{\sum^{N}_{j=1}\exp(\ell_j)}, \!\!\!\quad
\sum^{N}_{i=1} \mathrm{g}_i = 1, \quad \!\!\!\!\!
\mathbf{O}_{T_d} = \sum^{N}_{i=1} \mathrm{g}_i \mathbf{O}_i
\label{eq:MOE_coef}
\end{align}

\noindent \ra{where $N$ is the total number of experts, $\ell_i$ denotes the raw router score for expert $i$, $O_i$ is the output of expert $i$, and $\mathbf{O}_{T_d}$ represents the final output of the dense router.}

\ra{As mentioned earlier, the routing operator can also be implemented sparsely. In this setting, only a subset of experts $T$, containing the $K$ experts with the highest probability scores produced by the gating function, is evaluated \cite{fedus2022switch}. Let \(T=\mathrm{TopK}(\bm{\ell})\) denote the set of selected experts and let $O_{T_s}$ denote the final output. The renormalized sparse router is then defined as:}

\begin{align}
    \mathrm{\tilde{g}}_i= \frac{\exp(\ell_i)}{\sum_{j\in T}\exp(\ell_j)}, \quad \!\!\!\!
    \sum^{N}_{i=1} \mathrm{\tilde{g}}_i = 1, \quad  \!\!\!\!\!\!
    \mathbf{O}_{T_s}=\sum_{i\in T}\mathrm{\tilde{g}}_i \mathbf{O}_i
\end{align}

\ra{An MoE framework is well-suited for handling FDIAs, which are characterized by diverse shifts in operating regimes, because it enables different submodels to specialize and adaptively combine their outputs to address various cyberattacks. In addition, the MoE architecture supports joint detection and localization as a multi-task learning framework through input-conditioned gating \cite{kecceci2025federated, ma2018modeling}.}

\subsection{Encoder-only Transformer}

\ra{Typically, an MoE router is designed to be lightweight to reduce model complexity and improve inference speed \cite{lepikhin2020gshard}. However, this design can make the gating function highly sensitive to the input data. When the input is difficult to generalize using simple data-driven methods, the resulting probabilities may fail to reflect the optimal expert distribution. Therefore, preprocessing the input to generate richer feature representations can significantly improve routing performance. One effective approach is to provide globally informed representations using the Attention mechanism. Several studies have shown that combining Attention mechanisms with MoE architectures can achieve strong performance \cite{qiu2024layerwise}.}

\ra{First introduced in \cite{vaswani2017attention}, the Transformer architecture has demonstrated high performance in representation-learning and sequence modeling tasks \cite{khan2022transformers}. The model uses the attention mechanism to update each token representation by selectively aggregating information from other tokens in the sequence. These contextualized representations support tasks such as autoregressive generation, retrieval of relevant information, and feature embedding \cite{shen2016attention}.}

\ra{For FDIA detection and localization, contextualized embeddings produced by the attention block are particularly useful because they help identify inconsistencies between buses \cite{pan2026transformer}. Although the original Transformer architecture also includes a cross-attention layer in the decoder, its primary role is autoregressive decoding of encoder outputs, which is unnecessary for fixed-size classification tasks \cite{vaswani2017attention}.}

\ra{The scaled dot-product self-attention mechanism takes as input a sequence $\mathbf{X} \in \mathbb{R}^{n \times d_{model}}$ containing tokens (i.e., buses in this study). It consists of three main components obtained by linearly projecting $\mathbf{X}$ using learned weight matrices $\mathbf{W}_Q$, $\mathbf{W}_K$, and $\mathbf{W}_V$ with dimensions $d_{model}\times d_k$ for queries and keys, and $d_{model}\times d_v$ for values:}
\ra{\begin{itemize}
    \item Query ($\mathbf{Q}=\mathbf{X}\bm{W}_Q$), where $q_i$ represents the matching request of token $i$.
    \item Key ($\mathbf{K}=\mathbf{X}\bm{W}_K$), where $k_j$ represents the descriptor provided by token $j$ for matching.
    \item Value ($\mathbf{V}=\mathbf{X}\bm{W}_V$), where $v_j$ contains the information aggregated from token $j$.
\end{itemize}}

\ra{First, the similarity scoring (matching) step is performed between queries and keys, where each token’s query vector ($q_i$) is compared with the key vectors of other tokens ($k_j$) to determine which tokens are most relevant for updating the representation of token $i$:}

\begin{equation}
    \mathrm{Score}(\mathbf{Q},\mathbf{K}) = \frac{\mathbf{QK}^T}{\sqrt{d_k}}
    \label{eq:attn_similarity_score}
\end{equation}

\ra{The similarity scores are then normalized across keys for each query using a softmax activation function. The resulting coefficient matrix is used to compute a weighted sum of the value vectors, producing contextualized representations for each token:}

\begin{equation}
    Attention(\mathbf{Q,K,V}) = \mathrm{softmax}(\mathrm{Score}(\mathbf{Q,K}))\mathbf{V}
\label{eq:attn}
\end{equation}

\ra{One limitation of the scaled dot-product attention mechanism is that each token receives only a single attention distribution over all other tokens. Multi-head attention addresses this limitation by computing attention $h$ times using different learned projections and concatenating the resulting outputs into a single matrix. Under this formulation, the dimensions satisfy $d_k=d_v=d_{model}/h$.} \ra{Because these attention computations are independent, they can be executed in parallel, significantly improving computational efficiency \cite{fu2024transformer}. This process is illustrated in Figure \ref{fig:MultiHeadAttention} and can be formally expressed as:}

\begin{equation}
\!\! MultiHead(\mathbf{X}) = Concat(head_1, ..., head_h)\mathbf{W}^O
\label{eq:attn_mha}
\end{equation}

\begin{equation}
\mathrm{head}_i = \mathrm{Attention}(\mathbf{XW}_i^Q,\mathbf{XW}_i^K,\mathbf{XW}_i^V)
\end{equation}

{\ra{\noindent where $\mathbf{W}^O \in \mathbb{R}^{h\times d_v \times d_{model}}$.}

\ra{Although the original Transformer architecture includes both Encoder and Decoder components, this work focuses on the Encoder because it employs the self-attention mechanism to compute contextualized embeddings while requiring lower computational cost \cite{vaswani2017attention}. The Encoder consists of a stack of $N$ identical layers, each containing two sublayers: a multi-head attention mechanism and a position-wise fully connected feed-forward network. Each sublayer is connected through residual connections, and each layer, except the first, takes the output of the previous layer as input.}
\ra{In addition, the Encoder includes embedding and positional encoding (PE) modules to incorporate token order information. The embeddings and all sublayers share the same dimensionality $d_{model}$ \cite{vaswani2017attention}.}

\begin{figure}
    \centering
    \includegraphics[width=0.45\textwidth]{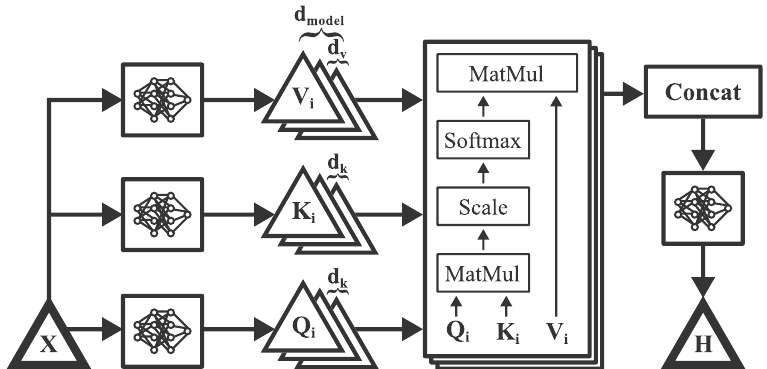}
    \caption{Multi-Head scaled dot-product attention mechanism. The input $X\in \mathbb{R}^{n\times d_{model}}$ is linearly projected to obtain queries (Q), keys (K), and values (V), which are reshaped into $h$ heads with per-head dimensions $d_k$ and $d_v$ (typically, $d_k=d_v$), summing to $d_{model}$. Then, each matrix is sent to one of the $h$ attention mechanisms depicted as a large rectangle with multiple functions inside. The attention mechanism illustrates the Eq.~\ref{eq:attn} and yields a contextualized output found in different subspaces. The head outputs then follow the Eq.~\ref{eq:attn_mha} through linear projection and concatenation into the context-aware final representation $H \in \mathbb{R}^{n\times d_{model}}$.}
    \label{fig:MultiHeadAttention}
\end{figure}

\subsection{Enabling topology awareness}

\ra{The traditional Transformer architecture does not naturally support graph-structured inputs. However, incorporating power network topology is essential for practical FDIA detection and localization. To adapt the sequence-based Transformer to graph-structured data, we integrate the Laplacian operator into the PE component and modify the attention formulation to explicitly incorporate graph connectivity through a scaled diffusion bias~\eqref{eq:attn_diff_bias}.}

\ra{To construct the Laplacian positional encoding (LapPE) block, let $\mathbf{L}_{sym}=\mathbf{I}-\boldsymbol{\Delta}^{-1/2}\mathbf{A}_w\boldsymbol{\Delta}^{-1/2}$ denote the symmetric normalized Laplacian, where $\mathbf{I}$ is the identity matrix and $\boldsymbol{\Delta}=\mathrm{diag}(d_1,d_2,\ldots,d_n)$ is the degree matrix with entries $d_i=\sum_{j=1}^{n}(\mathbf{A}_w)_{ij}$. Let $(\lambda_\ell,\bm{u}_\ell)$ represent the eigenvalue-eigenvector pairs of $\mathbf{L}{sym}$ ordered as $0=\lambda_1\leq \lambda_2\leq \cdots \leq \lambda_n$. The LapPE matrix is then defined as:}

\begin{equation}
    \mathbf{P} = [\bm{u}_2,\ldots,\bm{u}_{k_{\mathrm{PE}}+1}] \in \mathbb{R}^{n\times k_{\mathrm{PE}}},
\end{equation}

\noindent \ra{where the trivial constant eigenvector is excluded, and $\bm{u}_{k_{\mathrm{PE}}}$ denotes the $k_{\mathrm{PE}+1}$-th eigenvector, which captures global structural patterns across nodes \cite{shuman2013emerging}.}

\ra{The LapPE component injects graph-structural information into node representations, providing structural context for the $\mathbf{Q}$, $\mathbf{K}$, and $\mathbf{V}$ elements. However, using LapPE alone does not explicitly enforce graph-connectivity awareness in the attention similarity scores. Instead, the model must learn to use the LapPE features to capture topology, which is not guaranteed and can be ineffective during early training stages. To address this limitation, we introduce a graph-based diffusion matrix $\hat{D}$ as an additive bias in the similarity score formulation~\eqref{eq:attn} \cite{atwood2016diffusion}. The matrix is constructed from a $q$-step row-stochastic random-walk transition matrix.}

\begin{equation}
    \mathbf{D}^{q}=(\mathbf{T}_{rw}^{-1} \mathbf{A}_w)^{q}
    \label{eq:rw}
\end{equation}

\noindent \ra{where $\mathbf{T}_{rw}$ denotes the probability of reaching bus $j$ from bus $i$ in $q=4$ random-walk steps. Note that the matrix $\mathbf{D}$ is not directly suitable for addition to the pre-softmax similarity logits. Therefore, it is transformed into a multiplicative prior by applying a logarithmic transformation to each entry and standardizing the result to obtain $\mathbf{\hat{D}}$.}
\ra{Because adding a fixed bias may blur subtle discrepancies and potentially hide informative contrasts between nodes $i$ and $j$, we introduce a learnable scaling parameter $\gamma$ that allows the model to adaptively control the contribution of the diffusion bias. The attention formulation then becomes:}

\begin{equation}
    \mathrm{GraphAttn} = \mathrm{softmax}(\mathrm{Score}(\mathbf{Q,K})+\gamma \hat{\mathbf{D}})\mathbf{V}
    \label{eq:attn_diff_bias}
\end{equation}

\subsection{Feed-Forward Network as a router}

\ra{Due to router design constraints, the routing architecture must balance (i) routing quality, (ii) training stability, and (iii) computational complexity \cite{lepikhin2020gshard}, which is particularly important in online FDIA handling. While advanced methods such as attention mechanisms can enrich feature representations before routing, simpler approaches such as linear routers often lack sufficient expressiveness. Therefore, we adopt an MLP-based gating function, which provides expressive routing decisions with low latency \cite{harvey2025optimizing}.}

\ra{Widely used in MoE architectures \cite{shazeer2017outrageously, fedus2022switch, jiang2024mixtral}, feed-forward networks (FFNs) are well suited for three main reasons: (i) they improve routing quality by learning nonlinear mappings from tokens (electrical nodes) to expert contribution coefficients while remaining lightweight \cite{lo2025closer}; (ii) they help mitigate routing-related training instabilities because stabilization techniques, such as auxiliary differentiable load-balancing losses, are inexpensive and easy to integrate \cite{zhou2022mixture}; and (iii) they maintain low routing complexity by following the Transformer position-wise computation pattern, enabling efficient parallelization with a computational cost close to linear routers \cite{riquelme2021scaling}.}

\ra{Formally, a standard two-layer FFN $g(x)$ connected through an activation function $\sigma$ is defined as:}

\begin{equation}
g(x)= \mathrm{softmax}(\sigma(x\mathbf{W}_1+b_1)\mathbf{W}_2+b_2)
\label{eq:ffn_standard}
\end{equation}

\subsection{Architecture of the Proposed method}

\ra{The four modules are integrated into a unified framework consisting of: (i) residually connected LapPE for topology-aware electrical node embeddings, (ii) a scaled diffusion bias-enhanced attention mechanism that generates global contextualized representations of each bus for the localized MoE, (iii) a dense MoE with convolution-based experts that analyze each node through its neighborhood, and (iv) an MLP-based classifier that computes node-level FDIA binary flags, which are then aggregated for graph-level cyberattack detection. The overall architecture is illustrated in Figure \ref{fig:EOTConvMoE}.}

\ra{The model receives the node-feature matrix $\mathbf{X}^0 \in \mathbb{R}^{n \times 2}$, which contains the active ($P$) and reactive ($Q$) power injection measurements at each bus, along with a nonnegative weighted adjacency matrix $\mathbf{W} \in \mathbb{R}^{n \times n}$. Voltage magnitude $V$ and phase angle $\theta$ are intentionally excluded to avoid delays associated with PSSE output generation and to reflect practical operating conditions, since PMUs are not typically installed at every bus.}

\ra{First, the raw features $\mathbf{X}^{0}$ are projected to a $d_{model}$-dimensional space and combined with the fixed LapPE structural encoding precomputed from the graph and padded with zeros to match the projected shape of the features. A learnable scaling parameter balances the contribution of the structural information and the metering data.}

\begin{equation}
  \mathbf{H}^{0} = \mathrm{LinearNorm} (\mathbf{X}^{0}\mathbf{W}_x
    + \alpha \mathbf{P})
\end{equation}

\ra{The resulting matrix $\mathbf{H}^{0}$ is processed by diffusion-biased multi-head self-attention to generate contextualized representations, followed by a residual connection and normalization layer. The embedded representations $\mathbf{H} \in \mathbb{R}^{n \times d_{model}}$ are then used by the FFN router to produce expert weights $\mathrm{g}_1$ and $\mathrm{g}_2$. In parallel, $\mathbf{H}$ and $\mathbf{W}$ are passed to each expert to generate graph feature representations $\mathbf{O}_1 \in \mathbb{R}^{n \times d_{model}}$ and $\mathbf{O}_2 \in \mathbb{R}^{n \times d_{model}}$.}

\ra{We select graph-convolutional operators as experts because their topology-aware aggregation and weight-sharing mechanisms are well suited for modeling power-system measurements \cite{defferrard2016convolutional}. Specifically, we use (i) a CGCN expert, which applies localized spectral filtering through Chebyshev polynomial approximations to capture smooth spatial patterns \cite{boyaci2022cyberattack}, and (ii) an ARMAConv expert, which uses recursive graph filtering to model longer-range and noisier topological interactions. These complementary filtering behaviors are effective for FDIA detection, where anomalies may appear as both gradual and abrupt spatial variations \cite{boyaci2021joint}. The number of experts is kept small to limit computational complexity while preserving expert specialization.}

\ra{The graph feature representations $\mathbf{O}_1$ and $\mathbf{O}_2$ are weighted by the corresponding expert weights $\mathrm{g}_1$ and $\mathrm{g}_2$ and summed to produce the MoE output $\mathbf{O}_T \in \mathbb{R}^{n \times d{model}}$. This output is then added to and normalized with $\mathbf{H}$. The resulting representation $\mathbf{O}_{\mathrm{raw}}$ can be passed to another attention-enhanced convolution-based MoE layer for deeper processing or linearly projected to obtain cyberattack probabilities. A sigmoid function with a threshold of $0.5$ converts these probabilities into binary node-level attack labels $\mathbf{O}_L \in \mathbb{R}^{n}$. Finally, the graph-level FDIA detection output $\mathbf{O}_D \in \mathbb{R}^{1}$ is determined by taking the maximum value across all entries in $\mathbf{O}_L$.}

\begin{figure}[t]
    \centering
    \includegraphics[width=0.5\textwidth]{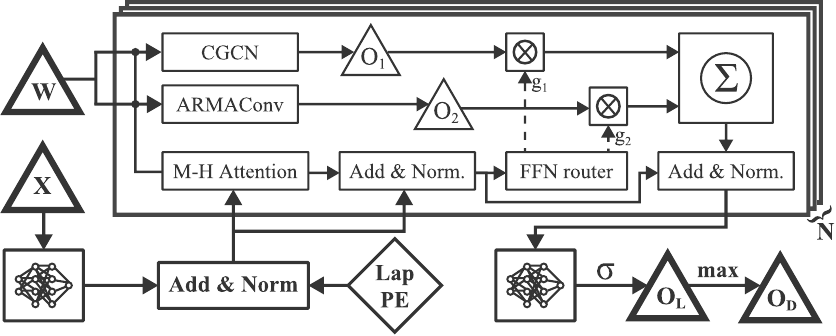}
    \caption{The proposed joint detector and localizer EOTConvMoE framework with dense router configuration in N attention-enhanced convolution-based MoE layers. The model takes the weighted adjacency matrix W and the raw power injections $P, Q$, compressed into $X$, as input and outputs both node-level binary attack flags and graph-level FDIA status. Note that the "Add\&Normalize" block for the linearly projected input and the LapPE component includes a learnable scalar to balance the contribution of the topological information.}
    \label{fig:EOTConvMoE}
\end{figure}

\section{Experimental setup and discussions}
\label{Experiments}

\ra{In this section, we describe the experimental setup, including the dataset generation process, and evaluate the performance of EOTConvMoE against common deep-learning benchmarks. For reproducibility and transparency, the codebase is publicly available at \url{https://github.com/netiRussell/FDIA-Transformer-MoE}.
}

\subsection{Dataset generation}

\ra{Due to security restrictions, real-world FDIA datasets are not publicly available. To address this limitation, we generate a synthetic dataset using the historical NYISO load profile from July 2021 and publicly available benchmark power-system networks to simulate false data injection attacks.}

\ra{We follow the benign dataset generation procedure in \cite{boyaci2021joint} to create intact 1-minute snapshots representing normal operating conditions for the IEEE 118- and 300-bus test systems. The resulting dataset contains $\frac{60\mathrm{min}}{1\mathrm{hour}} \times \frac{24\mathrm{hours}}{1\mathrm{day}} \times 31\mathrm{days} = 44{,}640$ samples.}

\ra{The final dataset is generated by modifying subsets of the original load measurements to match the FDIA patterns listed in Table \ref{tab:implemented_fdias}, while maintaining the attack distribution shown in Table \ref{tab:dataset_split}. First, a time step is randomly selected without replacement from a uniformly shuffled set of network states. Next, a bus is randomly chosen as the attack center, and a Breadth-First Search (BFS) algorithm identifies neighboring buses within $r$ hops, where $3 \le r \le 4$ for the IEEE 118-bus system and $6 \le r \le 8$ for the IEEE 300-bus system. Generator buses and zero-injection buses are excluded to preserve realistic operating conditions. The active and reactive power measurements of the selected buses are then modified according to the FDIA type, and the resulting matrices are stored. The $(\mathcal{A}_s)$ attack is sampled independently for each bus, while the $(\mathcal{A}_r)$ attack assumes $\tau=4$. The final dataset contains 34,560 samples, normalized using the unattacked samples from the training subset.}

\begin{table}[h!]
\centering
\caption{Implemented FDIAs}
\label{tab:implemented_fdias}
\resizebox{\columnwidth}{!}{
\begin{tabular}{|c|c|c|}
\hline
\textbf{FDIA type} & \textbf{Formulation} & \textbf{Reference} \\
\hline
\textbf{optimization-based} $(\mathcal{A}_o)$ & Eq.~(5) in \cite{boyaci2021graph}& \cite{boyaci2021joint} \\
\hline
\textbf{distribution-based} $(\mathcal{A}_d)$ & $\mathbf{z}_a(t) \sim \mathcal{N}\left( \mu(\mathbf{z}_o), \sigma^2(\mathbf{z}_o) \right)$ & \cite{chu2018unobservable} \\
\hline
\textbf{data scale} $(\mathcal{A}_s)$ & $\mathbf{z}_a(t) = \mathcal{U}(0.9, 1.1) \cdot \mathbf{z}_o(t).$ & \cite{paudel2024evaluation} \\
\hline
\textbf{data replay} $(\mathcal{A}_r)$ & $\mathbf{z}_a(t) = \mathbf{z}_o(t - \tau)$ & 
\cite{kecceci2025federated}\\
\hline
\end{tabular}
}
\end{table}

\begin{table}[h!]
\centering
\caption{Number of data points in each subset}
\label{tab:dataset_split}
\resizebox{\columnwidth}{!}{
\begin{tabular}{| c | c c c c c | c |} 
 \hline
 Subset & $A_{o}$ & $A_{d}$ & $A_{s}$ & $A_{r}$ & Non-attacked & Total \\
 \hline\hline
 Training & 5760 & 5760 & 0 & 0 & 11520 & 23040 \\ 
 \hline
 Validation & 1440 & 1440 & 0 & 0 & 2880 & 5760 \\ 
 \hline
 Testing & 720 & 720 & 720 & 720 & 2880 & 5760 \\ 
 \hline
\end{tabular}
}
\end{table}

\subsection{Training pipeline}

\ra{To evaluate the proposed model fairly, we compare it against several ML-based benchmarks, including traditional, temporal, convolutional, and attention-based architectures. Specifically, we implement MLP, LSTM, CGCN, ARMAConv, and ACEOT models following \cite{boyaci2022cyberattack, boyaci2021joint, abdulin2026attention}.}
\ra{To ensure a fair comparison, hyperparameters for each model are optimized under a fixed search budget. Each algorithm is tuned over 250 Optuna trials using the TPESampler and HyperbandPruner. The explored model-specific hyperparameters and selected values are summarized in Table \ref{tab:model_hyperparameters}, while the general training hyperparameters, including learning rate (LR), dropout, positive class weight, and LR warmup steps, are listed in Table \ref{tab:general_hyperparameters}.}

\ra{The selected hyperparameters are then used during training, which is formulated as a multi-label supervised learning problem. All models are trained using the AdamW optimizer and binary cross-entropy (BCE) loss for up to 256 epochs with mini-batches of 256 samples. Following \cite{vaswani2017attention}, the learning rate is scheduled with an initial warmup phase, where each step corresponds to one mini-batch, followed by gradual decay to zero using a cosine annealing schedule.}
\ra{To reduce overfitting, early stopping is enabled after the first 30 epochs. Training terminates when the validation loss fails to decrease by more than $10^{-4}$ over 16 consecutive epochs.}

\ra{While single-operator models can be trained using only BCE loss, MoE-based frameworks often require additional objectives to promote stable convergence and balanced expert utilization. In this work, we use the auxiliary loss $\mathcal{L}_{ent}$ to maximize the entropy of the batch-averaged routing distribution and encourage balanced expert usage \cite{chen2025mixture}. We also apply the supplemental loss $\mathcal{L}_z$ to penalize excessively large router logits and improve numerical stability during training \cite{zoph2022st}. The contribution $\mathrm{w_{loss}}$ of each loss term is learned during training, and the weighted losses accumulated across MoE layers are added to the classification loss.}

\begin{equation}
    \mathcal{L}_{ent} = w_{\mathrm{ent}}\sum_{e=1}^{E} \bar{p}_e \log\!\left(\bar{p}_e + \varepsilon\right)
\end{equation}

\noindent \ra{where $\bar{p}_e$ denotes the average router probability for expert $e$ across all graphs in the batch. When minimized with positive $w_{\mathrm{ent}}$, this term supports balanced router distribution}

\begin{equation}
    \mathcal{L}_{z}=w_z\frac{1}{B}\sum_{b=1}^{B}[\log (\sum_{e=1}^{E}\exp\big(z^{(b)}_{e}\big))]^2
\end{equation}

\noindent \ra{where $B$ denotes the batch size and $z^{(b)}_{e}$ represents the router logit assigned to expert $e$ for graph $b$.}

\ra{Empirical evidence suggests that pretraining experts can significantly accelerate the convergence of MoE-based models \cite{zhu2024moe}. Therefore, we modify the original ARMAConv and CGCN frameworks to align with the EOTConvMoE pipeline using a single operator. The resulting architectures are similar to Figure \ref{fig:EOTConvMoE}, but exclude the MoE-specific components, including the FFN, contribution weighting, expert output aggregation, and the associated ``Add\&Norm'' blocks. The tuned $\mathrm{ARMAConv}_{\mathrm{exp}}$ and $\mathrm{CGCN}_{\mathrm{exp}}$ models are then trained independently, and their learned weights are used to initialize the experts in EOTConvMoE during both hyperparameter tuning and model training.}

\begin{table}[t]
\centering
\caption{Model-specific hyperparameter tuning results}
\label{tab:model_hyperparameters}
\footnotesize
\setlength{\tabcolsep}{3pt}
\renewcommand{\arraystretch}{1.05}

\begin{tabularx}{\columnwidth}{@{} l l X c c @{}}
\hline
\textbf{Model} & \textbf{Param} & \textbf{Options} & \textbf{118} & \textbf{300} \\
\hline
MLP      & Layers      & \{2,3,4\}              & 4 & 3 \\
         & Units       & \{32,64,128\}          & 128 & 128 \\
LSTM     & Layers      & \{2,3,4\}              & 3 & 3 \\
         & Units       & \{32,64,128\}          & 128 & 128 \\
CGCN     & Layers      & \{2,3,4\}              & 4 & 4 \\
         & K           & \{3,4,5,7\}            & 7 & 7 \\
         & Units       & \{32,64,128\}          & 128 & 128 \\
ARMAConv & Layers      & \{2,3,4\}              & 3 & 3 \\
         & Stacks      & \{3,4,5\}              & 5 & 4 \\
         & Iterations  & \{4,5,6\}              & 6 & 5 \\
         & Units       & \{32,64,128\}          & 128 & 128 \\
ACEOT    & ARMA layers & \{2,3,4\}              & 3 & 4 \\
         & Stacks      & \{3,4,5\}              & 3 & 5 \\
         & Iterations  & \{4,5,6\}              & 6 & 5 \\
         & ARMA Units  & \{32,64,128\}          & 96 & 96 \\
         & T-heads     & \{2,4,6,8\}            & 2 & 2 \\
         & T-layers    & \{1,2,3\}              & 3 & 3 \\
         & $d_{\text{model}}$ & \{32,64,128\}   & 32 & 64 \\
         & $d_{\mathrm{FFN}}$ & \{64,128,256,512\}  & 64 & 128 \\
$\mathrm{CGCN}_{\mathrm{exp}}$ 
         & Layers      & \{2,3,4\}              & 4 & 4 \\
         & K           & \{3,4,5,7\}            & 7 & 7 \\
         & Units       & \{32,64,128\}          & 128 & 128 \\
$\mathrm{ARMAConv}_{\mathrm{exp}}$ 
         & Layers      & \{2,3,4\}              & 4 & 3 \\
         & Stacks      & \{3,4,5\}              & 4 & 4 \\
         & Iterations  & \{4,5,6\}              & 5 & 4 \\
         & Units       & \{32,64,128\}          & 128 & 128 \\
$\mathrm{EOTConvMoE_d}$   & T-heads     & \{2,4,6,8\}              & 4 & 4 \\
             & MoE-layers    & \{1,2,3,4\}              & 2 & 2 \\
             & $d_{\text{model}}$ & \{64,128,256\}   & 128 & 128 \\
             & $d_{\mathrm{FFN}}$ & \{64,128,256,512\}  & 256 & 256 \\
$\mathrm{EOTConvMoE_s}$   & T-heads     & \{2,4,6,8\}              & 4 & 4 \\
             & MoE-layers    & \{1,2,3,4\}              & 1 & 2 \\
             & $d_{\text{model}}$ & \{64,128,256\}   & 128 & 128 \\
             & $d_{\mathrm{FFN}}$ & \{64,128,256,512\}  & 512 & 256 \\
\hline
\end{tabularx}
\end{table}

\begin{table*}[t]
\centering
\caption{General hyperparameter tuning results}
\label{tab:general_hyperparameters}
\footnotesize
\setlength{\tabcolsep}{2pt}
\begin{tabularx}{\textwidth}{@{} l *{8}{>{\centering\arraybackslash}X} @{}}
\hline
\textbf{Model}
& \multicolumn{2}{c}{\textbf{LR}} & \multicolumn{2}{c}{\textbf{Dropout}}
& \multicolumn{2}{c}{\textbf{Positive weight}} 
& \multicolumn{2}{c}{\textbf{\# of Warmup steps}} \\
& \textbf{118} & \textbf{300}
& \textbf{118} & \textbf{300}
& \textbf{118} & \textbf{300}
& \textbf{118} & \textbf{300} \\
\hline
MLP         & $2.38\times10^{-3}$ & $2.89\times10^{-3}$ & 0.100 & 0.103 & 4.92 & 3.33 & 500 & 300 \\
LSTM        & $9.93\times10^{-3}$ & $8.14\times10^{-3}$ & 0.101 & 0.123 & 5.28 & 1.69 & 600 & 300 \\
CGCN        & $3.04\times10^{-3}$ & $1.00\times10^{-2}$ & 0.103 & 0.115 & 2.23 & 3.15 & 100 & 200 \\
ARMAConv    & $8.50\times10^{-4}$ & $1.01\times10^{-3}$ & 0.102 & 0.109 & 3.47 & 2.26 & 200 & 600 \\
ACEOT       & $3.43\times10^{-3}$ & $1.53\times10^{-3}$ & 0.119 & 0.181 & 2.23 & 2.87 & 100 & 500 \\
$\mathrm{CGCN}_{\mathrm{exp}}$  & $1.18\times10^{-2}$ & $1.03\times10^{-2}$ & 0.109 & 0.102 & 3.99 & 3.70 & 200 & 100 \\
$\mathrm{ARMAConv}_{\mathrm{exp}}$  & $3.38\times10^{-3}$ & $2.91\times10^{-3}$ & 0.123 & 0.108 & 4.28 & 4.72 & 200 & 400 \\
$\mathrm{EOTCMoE_d}$  & $4.22\times10^{-4}$ & $5.73\times10^{-4}$ & 0.104 & 0.101 & 6.79 & 4.57 & 500 & 300 \\
$\mathrm{EOTCMoE_s}$  & $2.17\times10^{-3}$ & $1.50\times10^{-3}$ & 0.129 & 0.111 & 4.11 & 4.66 & 300 & 300 \\ 
\hline
\end{tabularx}
\end{table*}

\subsection{Evaluation techniques and metrics}

\ra{Given the classification nature of the problem, we employ confusion matrix-based metrics derived from True Positive (TP), False Positive (FP), True Negative (TN), and False Negative (FN) predictions. Specifically, we evaluate the detection rate (DR), false alarm rate (FA), and F1-score, defined as $DR = \frac{\mathrm{TP}}{\mathrm{TP} + \mathrm{FN}}$, $FA = \frac{\mathrm{FP}}{\mathrm{FP} + \mathrm{TN}}$, and $F1 = \frac{2\mathrm{TP}}{2\mathrm{TP} + \mathrm{FN} + \mathrm{FP}}$, respectively.}
\ra{In addition, we employ mean Average Precision (mAP), computed as the mean AP over all attacked samples. For each attacked sample, buses are ranked in descending order according to their predicted attack scores. The AP is defined as $\mathrm{AP} = \frac{1}{P} \sum_{k=1}^{n} \mathrm{Precision}(k)y_k$, where $\mathrm{Precision}(k)=\frac{1}{k}\sum_{i=1}^{k}y_i$. Here, $P$ denotes the number of attacked buses in the sample, and $y_k$ represents the ground-truth label of the bus ranked at position $k$.}

\ra{To avoid division-by-zero issues in benign samples, we assign $DR=1$, $FA=0$, and $F1=1$ when all buses are correctly classified; otherwise, we set $DR=0$, $FA=1$, and $F1=0$. Benign samples are excluded from the mAP computation because they correspond to $P=0$, which makes AP undefined. Confusion matrix-based metrics are used due to their interpretability and direct relevance to classification performance, while mAP is included to evaluate ranking quality independently of a specific decision threshold.}

\ra{This work evaluates the performance of the proposed EOTConvMoE model through several comparative experiments. Specifically, we investigate the impact of the routing mechanism by evaluating the densely routed $\mathrm{EOTConvMoE_d}$ and sparsely routed $\mathrm{EOTConvMoE_s}$ variants. We further analyze attack detection performance using the F1-score, false alarm rate (FA), and detection rate (DR), and assess attack localization capability using the sample-wise F1-score and mean Average Precision (mAP). Finally, computational efficiency is evaluated based on the wall-clock inference time required to generate model outputs.}

\subsection{Numerical results and Analysis}

\ra{The experiments are conducted on a workstation equipped with an Intel(R) Xeon(R) w7-2595X CPU operating at 2.81 GHz, an NVIDIA RTX PRO 6000 Blackwell Max-Q GPU, and 128 GB of RAM. Our codebase is implemented in Python 3.13 using the \texttt{PyTorch} and \texttt{Pandapower} libraries.}

\ra{Table \ref{tab:detection_table} presents the performance of EOTConvMoE and the benchmark architectures. The results show that leveraging multiple experts can substantially improve detection performance compared to single-operator methods. Although the MLP model, which is the simplest architecture considered, achieves the lowest F1 score, temporal- and convolution-based operators demonstrate competitive performance. These findings suggest that combining attention-based mechanisms and convolutional experts can provide a strong attack detection capability.}

\ra{The proposed EOTConvMoE framework achieves the highest F1 score among all evaluated models. In particular, it outperforms the next best non-MoE baseline by 8.03\% and 6.41\% in the IEEE 118- and 300-bus test cases, respectively. Similarly, EOTConvMoE achieves the best DR values, exceeding the best non-MoE baseline by 13.19\% and 9.69\% in the IEEE 118- and 300-bus systems, respectively. In terms of the FA metric, the proposed framework achieves the best result for the IEEE 300-bus test case and remains only 0.07\% below the best score for the IEEE 118-bus test case. These results indicate that EOTConvMoE effectively balances attack detection capability with accurate identification of normal operating conditions.}

\ra{The observed performance improvements can be attributed to the dataflow pipeline embedded within the proposed architecture. Specifically, motivated by the benefits of combining encoder-only Transformer models with graph-filtering networks within a unified framework \cite{abdulin2026attention}, the input to the convolutional experts is preprocessed using topology-aware global information extracted by the attention mechanism. In addition, expert pre-training enables the convolutional operators to adapt more effectively to the transformed feature representations. The sparse routing configuration generally produces lower detection performance than the dense routing strategy because only a single filter type is activated at each layer.}

\begin{table}[t]
\centering
\caption{Detection results in DR, FA, and F1 percentages}
\label{tab:detection_table}
\footnotesize
\setlength{\tabcolsep}{3pt}
\renewcommand{\arraystretch}{1.05}
\begin{tabularx}{\columnwidth}{@{} l *{6}{>{\centering\arraybackslash}X} @{}}
\hline
\textbf{Model}
& \multicolumn{2}{c}{\textbf{F1 (\%)}} 
& \multicolumn{2}{c}{\textbf{DR (\%)}}
& \multicolumn{2}{c}{\textbf{FA (\%)}} \\

& \textbf{118} & \textbf{300}
& \textbf{118} & \textbf{300}
& \textbf{118} & \textbf{300} \\
\hline
MLP         & 50.03 & 34.92 & 40.42 & 25.10 & 21.15 & 18.66 \\
LSTM        & 84.98 & 85.22 & 75.35 & 74.66 & 1.98 & 0.56 \\
CGCN        & 83.55 & 83.17 & 71.77 & 71.27 & \bf{0.03} & \bf{0.14} \\
ARMAConv    & 85.88 & 85.53 & 75.42 & 75.50 & 0.21 & 1.04 \\
ACEOT       & 85.14 & 84.52 & 74.20 & 73.34 & 0.10 & 0.21 \\
$\mathrm{EOTConvMoE_d}$  & \bf{93.91} & \bf{91.94} & \bf{88.61} & \bf{85.19} & 0.10 & \bf{0.14} \\
$\mathrm{EOTConvMoE_s}$  & 93.13 & 90.92 & 87.22 & 83.51 & 0.11 & 0.19  \\ 
\hline
\end{tabularx}
\end{table} 

\ra{Table \ref{tab:localization_talbe} presents the attack localization results for the proposed and benchmark models. Similar to the detection analysis, EOTConvMoE and ACEOT achieve the highest F1 scores in both test cases, further indicating the effectiveness of combining attention-based and convolutional mechanisms within a unified framework. In addition, the observed strong correlation between the F1 and mAP scores suggests that the selected decision threshold produces a consistent ranking of model performance.}

\ra{Among the evaluated methods, CGCN, ACEOT, and EOTConvMoE achieve the highest localization performance in terms of mAP. In particular, EOTConvMoE reaches the highest mAP value of 84.82\% in the IEEE 118-bus test case and achieves the second-highest result in the IEEE 300-bus system, trailing the best model by only 0.07\%. Based on the F1 metric, the proposed framework achieves the best localization performance, outperforming the strongest non-MoE baseline, ACEOT, by 1.41\% and 1.01\% in the IEEE 118- and 300-bus test cases, respectively.}

\ra{The improved localization capability can be attributed to the complementary contributions of polynomial and rational graph filters within the expert ensemble. This interpretation is further supported by the performance similarity between the sparsely routed variant of EOTConvMoE and ACEOT, indicating the reduced flexibility associated with single-operator routing, even when the selected operator matches the local data characteristics. Convolution-driven architectures also demonstrate competitive performance, with all graph convolution-based models achieving F1 scores above 80\% across both test systems. In contrast, the conventional LSTM-based temporal model, despite lacking explicit graph information, still achieves F1 scores of 68.04\% and 64.17\% in the IEEE 118- and 300-bus systems, respectively. These results suggest that incorporating temporal experts into an MoE framework may further improve model diversity and representation capability.}

\begin{table}[t]
\centering
\caption{Localization results in F1 and mAP percentages}
\label{tab:localization_talbe}
\footnotesize
\setlength{\tabcolsep}{3pt}
\renewcommand{\arraystretch}{1.05}
\begin{tabularx}{\columnwidth}{@{} l *{4}{>{\centering\arraybackslash}X} @{}}
\hline
\textbf{Model}
& \multicolumn{2}{c}{\textbf{F1 (\%)}} 
& \multicolumn{2}{c}{\textbf{mAP (\%)}} \\

& \textbf{118} & \textbf{300}
& \textbf{118} & \textbf{300} \\
\hline
MLP         & 44.92 & 42.35 & 36.93 & 30.99 \\
LSTM        & 68.04 & 64.17 & 57.11 & 55.13 \\
CGCN        & 81.39 & 80.92 & 82.20 & 81.18 \\
ARMAConv    & 82.41 & 82.15 & 76.92 & 72.74 \\
ACEOT       & 83.56 & 82.72 & 82.70 & \bf{81.83} \\
$\mathrm{EOTConvMoE_d}$  & \bf{84.97} &\bf{83.73} & \bf{84.82} & 81.76 \\
$\mathrm{EOTConvMoE_s}$  & 83.60 & 82.33 & 83.21 & 80.47 \\ 
\hline
\end{tabularx}
\end{table}

\ra{Finally, we evaluate the inference latency of the benchmark models. For the IEEE 118-bus test case, the average inference times required to detect and localize a single test instance are 3.69 ms for MLP, 3.26 ms for LSTM, 19.02 ms for CGCN, 8.25 ms for ARMAConv, 9.65 ms for ACEOT, 31.31 ms for $\mathrm{EOTConvMoE_d}$, and 21.22 ms for $\mathrm{EOTConvMoE_s}$. For the IEEE 300-bus test case, the corresponding inference times are 3.74 ms, 3.40 ms, 20.46 ms, 7.46 ms, 16.77 ms, 65.62 ms, and 35.75 ms, respectively.}

\ra{The results indicate that the increased architectural complexity of the proposed framework, while improving F1 and mAP performance, also increases computational cost. In particular, the choice of expert operators significantly affects inference efficiency. For example, in the IEEE 118-bus system, CGCN, the slowest single-operator benchmark model, requires only 2.20 ms less inference time than the sparsely routed EOTConvMoE configuration. Additionally, larger graph sizes increase the latency of attention-based architectures, as observed for both ACEOT and EOTConvMoE in the IEEE 300-bus test case.}
\ra{Despite the higher computational cost, EOTConvMoE demonstrates good scalability, as its detection and localization performance remains stable when transitioning from the IEEE 118-bus system to the IEEE 300-bus system. Moreover, the computational overhead of the attention and expert aggregation operations can be mitigated through GPU parallelization. An important advantage of the MoE-based framework is its flexibility in balancing inference latency and classification performance through different routing strategies and expert selections. In particular, the sparsely routed $\mathrm{EOTConvMoE_s}$ achieves approximately 1.5$\times$ and 1.8$\times$ faster inference latency than the densely routed variant for the IEEE 118- and 300-bus systems, respectively, while maintaining competitive F1 and DR performance.}

\section{Conclusion}
\label{Conclusion}

\ra{Due to the increasing complexity of cyber-physical power networks, modern smart grids remain highly vulnerable to cyberattacks. Although numerous approaches have been proposed for false data injection attack (FDIA) detection, the use of multiple complementary operators within a unified architecture for joint attack detection and localization remains insufficiently explored. Unlike existing MoE-based FDIA mitigation methods, this work introduces a topology-aware graph-filter MoE framework in which the experts consist of graph convolutional filters with complementary characteristics, specifically Chebyshev and ARMA graph operators. In addition, the routing mechanism is enhanced through input preprocessing based on Laplacian positional encoding and diffusion-biased self-attention. Extensive simulations and comparative experiments conducted on the IEEE 118- and 300-bus benchmark systems using real-world NYISO load profiles demonstrate the effectiveness of the proposed EOTConvMoE framework. The results show that EOTConvMoE maintains low false-alarm rates, achieves superior F1 and DR detection performance compared to existing methods, and attains the highest localization F1 scores in both test systems. Furthermore, the proposed architecture adaptively learns expert contribution weights, enabling stable performance across different attack scenarios. Future work may investigate the incorporation of temporal modeling and physics-informed training strategies to further improve robustness and generalization capability.}

\bibliographystyle{IEEEtran}
\bibliography{ref}
\end{document}